\documentclass[11pt,letterpaper,reqno]{amsart}
\usepackage[left=30mm,top=30mm,right=30mm,bottom=30mm]{geometry}
\usepackage{epsfig}
\usepackage{amssymb}
\usepackage{amsthm}
\usepackage{subcaption}
\usepackage{multicol}
\usepackage{times}
\usepackage{color}
\usepackage{enumerate}
\usepackage{graphicx}
\usepackage{nicefrac}
\usepackage{xcolor}
\usepackage{tikz}
\usepackage[all]{xy}
\usepackage{comment}
\usepackage{hyperref}
\hypersetup{
 colorlinks=true, 
 linkcolor=blue, 
 citecolor=blue, 
 filecolor=blue, 
 urlcolor=blue 
}
\usepackage[nameinlink,capitalize,noabbrev]{cleveref}

\newtheorem{theorem}{Theorem}
\newtheorem*{theorem*}{Theorem}

\newtheorem{lemma}[theorem]{Lemma}

\newtheorem*{claim*}{Claim}

{\theoremstyle{remark}
 \newtheorem{remark}[theorem]{Remark}}
 \newtheorem*{remark*}{Remark}
{\theoremstyle{definition}

}

\newcommand{\Z}{\mathbb Z}

\newcommand{\K}{\mathbf{k}}

\begin{document}

\title[Non-normalized Nash blowup fails to resolve singularities in dimension three]{Non-normalized Nash blowup fails to resolve singularities \\ in dimension three}

\author[{F. Castillo}]{Federico Castillo}
\address{Facultad de Matemáticas, Pontificia Universidad Católica de Chile, Santiago, Chile}  
\email{federico.castillo@uc.cl}

 \author[D. Duarte]{Daniel Duarte}
\address{Centro de Ciencias Matematicas, UNAM, Morelia, M\'exico}  \email{adduarte@matmor.unam.mx}

\author[M. Leyton-Alvarez]{Maximiliano Leyton-\'Alvarez} %
\address{Instituto de Matem\'aticas, Universidad de Talca, Talca, Chile} %
\email{mleyton@utalca.cl}

\author[A. Liendo]{Alvaro Liendo} %
\address{Instituto de Matem\'aticas, Universidad de Talca, Talca, Chile} %
\email{aliendo@utalca.cl}

\date{\today}

\thanks{
 \mbox{\hspace{11pt}}{\it Key words}: Resolution of singularities, Nash blowup, Toric varieties.\\
 \mbox{\hspace{11pt}} The first author was partially supported by Fondecyt project 1221133. The second author was partially supported by CONAHCYT project CF-2023-G-33 and PAPIIT grant IN117523. The third author was partially supported by Fondecyt project 1221535. The fourth author was partially supported by Fondecyt project 1240101.  Finally, the first, third, and fourth authors were partially supported by Fondecyt Exploraci\'on project 13250049. }

\begin{abstract}
In this paper we show that iterating (non-normalized) Nash blowups does not necessarily resolve the singularities of algebraic varieties of dimension three over fields of characteristic zero.
\end{abstract}

\maketitle

% \tableofcontents

\section*{Introduction}

Let $X \subseteq \K^n$ be an equidimensional algebraic variety of dimension $d$, where $\K$ is an algebraically closed field. Consider the Gauss map:
\begin{align}\label{eq:gauss}
\Phi\colon X \setminus \operatorname{Sing}(X) \to \operatorname{Grass}(d, n) \quad \text{defined by} \quad x \mapsto T_x X\,,
\end{align}
where $\operatorname{Grass}(d, n)$ denotes the Grassmannian of $d$-dimensional vector spaces in $\K^n$, and $T_x X$ is the tangent space to $X$ at $x$. 
Let $X^*$ be the Zariski closure of the graph of $\Phi$ and $\nu\colon X^* \to X$ be the composition of the inclusion $X^*\hookrightarrow X\times\operatorname{Grass}(d, n)$ and the projection onto the first coordinate. The morphism $\nu$ is a proper birational map that is an isomorphism over $X \setminus \operatorname{Sing}(X)$. The pair $(X^*,\nu)$ is called the \emph{Nash blowup} of $X$. The composition $\eta\circ\nu$, where $\eta$ is the normalization map, is called the \emph{normalized Nash blowup} of $X$.

It has been conjectured that iterating the Nash blowup or the normalized Nash blowup resolves the singularities of an algebraic variety \cite{Se54,Sp90,GS1}.
The starting point for studying these questions is the theorem, due to A. Nobile, stating that the Nash blowup of $X$ is an isomorphism if and only if $X$ is non-singular, in characteristic zero; in prime characteristic this statement is false \cite{Nob75}. In particular, the question on the resolution properties of the (non-normalized) Nash blowup is studied only in characteristic zero.

There are many partial results giving an affirmative answer to the conjectures for some families of varieties \cite{Nob75,Reb,GS1,GS2,Hi83,Sp90,GS3,EGZ,Ataetal,GrMi12,GoTe14,DuarSurf,DG,DJNB,DDR,DS}. In a sudden change of course, the answer was ultimately negative in dimensions four and higher \cite{CDLAL}.

In view of the results of \cite{CDLAL} in the case of characteristic zero fields, it remains to find out the resolution properties of (non-normalized) Nash blowups in dimensions two and three. In this paper we exhibit a counterexample to this question in dimension three. More precisely, we prove the following theorem.

\begin{theorem*}\label{main}
For any field of characteristic zero, there exists a non-normal affine algebraic variety $X$ of dimension $3$ such that the second iteration of the Nash blowup of $X$ contains an open affine subset isomorphic to $X$.
\end{theorem*}

We prove the theorem by providing an explicit affine non-normal toric variety $X$ satisfying the hypothesis of the theorem. The study of Nash blowups of toric varieties over algebraically closed fields of characteristic zero was initiated by G. Gonz\'alez-Sprinberg in the seminal paper \cite{GS1}. He gave a combinatorial description of the normalized Nash blowup of a normal toric variety. This result was later strengthened by removing the assumption of normality, by P.~Gonz\'alez and B.~Teissier~\cite{GoTe14}, and by D.~Grigoriev and P.~Milman~\cite{GrMi12} (see also~\cite{LJ-R}). This combinatorial description serves as our main tool in proving the theorem.

Even though the proof of the theorem presented here is fully verifiable without the use of computational software, it is worth emphasizing that the affine non-normal toric variety $X$ used in the proof was discovered through extensive computer-assisted experimentation using SageMath \cite{sagemath}. 

This experimental approach to studying Nash blowups of toric varieties has been successfully employed in previous works such as \cite{Ataetal,CDLAL}. In both of these papers, the authors implement the normalized Nash blowup, which applies only to normal toric varieties and can therefore be treated using tools from convex geometry such as cones and fans. In contrast, the present paper fundamentally addresses non-normal toric varieties, which must be handled using the more general language of affine semigroups. As a result, the computational treatment required to obtain the counterexample presented here becomes significantly more involved.

As a final remark, we would like to note that we have also conducted extensive computational research on the two open cases of the Nash blowup conjecture: the non-normalized Nash blowup of two-dimensional varieties, and the normalized Nash blowup of three-dimensional varieties. To date, no counterexamples have been found, and it seems unlikely that such a counterexample will be discovered using toric methods.

\subsection*{Acknowledgments}

This collaboration began at the \href{https://sites.google.com/view/agrega0/home}{AGREGA} workshop, which took place at Universidad de Talca in January 2024 where the second named author presented the Nash blowup conjecture. We extend our gratitude to the institution for its support and hospitality.

\section{The Nash blowup of an affine toric variety}
\label{sec:toric}

A toric variety $X$ is a variety endowed with a faithful regular action of an algebraic torus having an open orbit. Toric varieties admit a well-known combinatorial description that we recall now, for details, see \cite{fulton1993introduction,oda1983convex,Sturm,cox2011toric}. In contrast to some references, but in line with \cite{Sturm,cox2011toric}, we do not require toric varieties to be normal.

Let $S$ be an affine semigroup, i.e., a finitely generated semigroup with identity element that can be embedded in a free abelian group. In the sequel, without loss of generality, we assume that $S\subset M=\Z^d$ and that the group generated by $S$ is $M$. 

We say that $S$ is saturated if for any $u\in M$ such that $\lambda u\in S$ for some $\lambda\in\Z_{\geq1}$, we have $u\in S$. We say that the semigroup $S$ is pointed if $S\cap(-S)=\{0\}$. The Hilbert basis of a pointed affine semigroup is its unique minimal generating set. It corresponds to elements $h\in S$ that cannot be written as a sum $h=m+m'$ with $m,m'\neq 0$.

Given an affine semigroup $S$, we define the semigroup algebra
\[
\K[S]=\bigoplus_{u\in S}\K\cdot\chi^{u},\quad\mbox{with}\quad\chi^{0}=1,\mbox{ and }\chi^{u}\cdot\chi^{u'}=\chi^{u+u'},\ \forall u,u'\in S\,.
\]
The affine variety $X(S):=\operatorname{Spec}\K[S]$ is a toric variety. Moreover, the affine toric variety $X(S)$  is  normal if and only if $S$ is saturated \cite[Theorem~1.3.5]{cox2011toric}.

\medskip

As mentioned in the introduction, the Nash blowup of a toric variety over algebraically closed fields of characteristic zero fields has a combinatorial description. We now recall this description following the account given in \cite[Section~1.9.2]{Sp20}.

Let $X(S)$ be the affine toric variety given by the pointed semigroup $S\subset \Z^d$ and let $H=\{h_1,\dots,h_r\}$ be a generating set of $S$. For a collection of $d$ elements $\{h_{i_1},\dots,h_{i_d}\}\subset S$, we define the matrix $(h_{i_1} \cdots h_{i_d})$ whose columns are the vectors $h_{i_j}$. 

The affine charts of the Nash blowup of $X$ are indexed by the subsets $A=\{h_{i_1},\dots,h_{i_d}\}$ of $S$ such that $\operatorname{det}(h_{i_1}\cdots h_{i_d})\neq 0\,.$
Let $A$ be such a subset of $S$. Without loss of generality, up to reordering the indices, we may and will assume $A=\{h_1,\dots,h_d\}$. Now let $h\in A$. Again, up to reordering the indices, we may and will assume $h=h_1$. We let
\begin{align} \label{eq:Ga}
\mathcal{G}_A(h)=\Big\{g-h\mid g\in H\setminus A \mbox{ and } \operatorname{det}(g\,h_2\cdots h_d)\neq 0 \Big\}.
\end{align}
Finally, letting
$\mathcal{G}_A=A\cup \mathcal{G}_A(h_1)\cup\dots\cup \mathcal{G}_A(h_d)$, we let  $S_A$ be the semigroup generated in $M$ by $\mathcal{G}_A$. Now,  the collection of affine toric varieties $X(S_A)$, for all $A=\{h_{i_1},\dots,h_{i_d}\}$ with 
$$\operatorname{det}(h_{i_1}\cdots h_{i_d})\neq 0 \quad \mbox{and}\quad S_A\mbox{ pointed}, $$
provides a set of covering affine charts of the Nash blowup of $X(S)$.

In the next section we will use this description to provide the example stated in the theorem for algebraically closed fields. For non-algebraically closed fields, see \cref{rem: non algeb clos}.

\section{Proof of the Theorem}

As stated in the introduction, to prove the theorem we will exhibit a non-normal affine toric variety of dimension 3 fulfilling the conditions of the theorem. Let $M=\Z^3$, and consider the semigroup $S\subset M$ generated by the columns of the matrix

$$
B=\left[\begin{array}{rrrrrr}
1 & 0 & 0 & -2 & 1 & 2 \\
0 & 1 & 0 & -1 & -1 & -2 \\
0 & 0 & 1 & 2 & 1 & 1
\end{array}\right]\,.
$$
Let $X(S)$ be the toric variety defined by $S$. We denote by $h_i$ the vector corresponding to the $i$-th column of $B$, with $i \in \{1, \dots, 6\}$. Hence $H=\{h_1,\ldots,h_6\}$ is a generating set of $S$.

\begin{lemma}
The affine semigroup $S$ is pointed and non-saturated.
\end{lemma}

\begin{proof}
Let $L(x,y,z)=x+2y+3z$. 
Then $L(h_i)\geq2$ for $i\in\{2,3,4,5\}$ and $L(h_1)=L(h_6)=1$. 
This implies that $S$ is pointed. 
Moreover, let $u=(0,-1,1)\in M$. 
Then $3u=(0,-3,3)=(-2,-1,2)+(2,-2,1)\in S$. 
But $u\notin S$ since $L(u)=1$.
\end{proof}

Let $A=\{h_1,h_4,h_6\}$. Notice that $\operatorname{det}(h_1\,h_4\,h_6)=3$. We have $H\setminus A=\{h_2,h_3,h_5\}$. There are 9 relevant determinants to be computed that we exhibit below (see (\ref{eq:Ga})). At the beginning of each line we show the element of $A$ that is omitted.
\begin{align} \label{eq:dets}
\begin{array}{llll}
h_1:&\operatorname{det}(h_2\,h_4\,h_6)=6, & \operatorname{det}(h_3\,h_4\,h_6)=6, &
\operatorname{det}(h_5\,h_4\,h_6)=3,
\\
h_4:&\operatorname{det}(h_2\,h_1\,h_6)=-1, & \operatorname{det}(h_3\,h_1\,h_6)=-2, &
\operatorname{det}(h_5\,h_1\,h_6)=-1,\\
h_6:&\operatorname{det}(h_2\,h_1\,h_4)=-2, & \operatorname{det}(h_3\,h_1\,h_4)=-1, &
\operatorname{det}(h_5\,h_1\,h_4)=1,
\end{array}
\end{align}

We obtain that $\mathcal{G}_A$ is
\begin{align*}
\mathcal{G}_A=A&\cup\{h_2-h_1,h_3-h_1,h_5-h_1\}\\ &\cup \{h_2-h_4,h_3-h_4,h_5-h_4\}\\&\cup\{h_2-h_6,h_3-h_6,h_5-h_6\}.\\
\end{align*}
Recall that $S_A$ denotes the semigroup generated by $\mathcal{G}_A$. The following subset $H_1$ of $\mathcal{G}_A$ is a generating set for $S_A$ 
$$H_1=\Big\{h_2-h_1,h_5-h_1,h_5-h_4,h_1,h_2-h_6,h_4,h_2-h_4,h_6\Big\}\,.$$
Indeed, for all the other elements in $\mathcal{G}_A$ we have
\begin{align*}
h_5-h_6&=h_2-h_1, \\
h_3-h_6&=2(h_2-h_1),  \\
h_3-h_1&=(h_2-h_1)+(h_5-h_1), \\
h_3-h_4&=(h_2-h_4)+(h_5-h_1).
\end{align*}

Consider the following matrix whose columns are the vectors in $H_1$
$$
B_1=\left[\begin{array}{rrrrrrrr}
-1 & 0 & 3 & 1 & -2 & -2 & 2 & 2 \\
1 & -1 & 0 & 0 & 3 & -1 & 2 & -2 \\
0 & 1 & -1 & 0 & -1 & 2 & -2 & 1

\end{array}\right]\,.
$$

Let $X(S_A)$ be the toric variety defined by $S_A$. Denote as $g_i$ the vector corresponding to the $i$-th column of $B_1$, with $i \in \{1, \dots, 8\}$.

\begin{lemma}
The affine semigroup $S_A$ is pointed. 
\end{lemma}
\begin{proof}
For $L_1(x,y,z)=5x+8y+10z$ we have $L_1(g_i)>0$ for each $i$.
\end{proof}

Let $A_1=\{g_1,g_5,g_7\}\subset H_1$. Notice that $\det(g_1\,g_5\,g_7)=-2$. We have $H_1\setminus A_1=\{g_2,g_3,g_4,g_6,g_8\}$. Repeating the previous computations we obtain
\begin{align*}
\mathcal{G}_{A_1}=A_1 &\cup\{g_2-g_1,g_3-g_1,g_4-g_1,g_6-g_1,g_8-g_1\}\\
&\cup\{g_2-g_5,g_3-g_5,g_4-g_5,g_6-g_5,g_8-g_5\}\\
&\cup\{g_3-g_7,g_4-g_7,g_6-g_7,g_8-g_7\}\,.
\end{align*}

Let $T$ be the semigroup generated by $\mathcal{G}_{A_1}$. The following subset $H_2$ of $\mathcal{G}_{A_1}$ generates $T$
$$H_2=\Big\{g_2-g_1,g_3-g_1,g_4-g_1,g_5,g_6-g_1,g_6-g_7\Big\}\,.$$
Consider the following matrix whose columns are the vectors in $H_2$
$$
B_2=\left[\begin{array}{rrrrrr}
1 & 4 & 2 & -2 & -1 & -4\\
-2 & -1 & -1 & 3 & -2 & -3\\
1 & -1 & 0 & -1 & 2 & 4
\end{array}\right]\,.
$$

\medskip

To conclude, let $U\colon M\to M$ be the automorphism given by the unimodular matrix
$$U=\left[\begin{array}{rrr}
1 & 4 & 2 \\
-2 & -1 & -1 \\
1 & -1 & 0
\end{array}\right]\,.$$

A straightforward verification yields that $UB=B_2$ and so $U$ induces a bijection from $H$ to $H_2$. 

This yields that $S$ is isomorphic to $T$. Since $S$ is pointed, the same holds for $T$. Hence $X(T)$ is an affine chart of the Nash blowup of $X(S_A)$, which in turn is an affine chart of the Nash blowup of $X(S)$. By the isomorphism $S\cong T$ we conclude that $X(S)\cong X(T)$. This achieves the proof of the theorem in the case where the base field is algebraically closed. 

\begin{remark}\label{rem: non algeb clos}
The combinatorial description that we used in the previous proof is based on the fact that the Nash blowup of a toric variety coincides with the blowup of the logarithmic Jacobian ideal. This was proved over algebraically closed fields \cite[Proposition 60]{GoTe14}. Hence, we cannot apply the same method over other characteristic zero fields. However, there is another approach that allows us to show that the counterexample also holds in that case. The following discussion is based on \cite[Construction 4.4]{GrMi12}. In what follows we use the notation established on the course of the proof.

Let $X(S)\subset\K^6$ be the 3-dimensional toric variety corresponding to $S=\mathbb{Z}_{\geq0}(h_1,\ldots,h_6)$. Recall that $X(S)$ can also be obtained as the Zariski closure of the image of the monomial map
$$\phi:(\K^*)^3\to\K^6,\,\,\,x\mapsto (x^{h_1},\ldots,x^{h_6}).$$
The tangent space $T_{\phi(x)}X(S)$ is determined by the Jacobian matrix  $J=J(x^{h_1},\ldots,x^{h_6})$. Using the Pl\"ucker embedding $P:\operatorname{Grass}(3,6)\hookrightarrow\mathbb{P}_{\K}^{\binom{6}{3}-1}$, the Nash blowup of $X(S)$ is then obtained as the Zariski closure of 
$$\{(\phi(x),(\ldots:\Delta_{i_1i_2i_3}:\ldots))\mid x\in(\K^*)^3\}\subset X(S)\times\mathbb{P}_{\K}^{\binom{6}{3}-1},$$
where $\Delta_{i_1i_2i_3}$ denotes the maximal minor defined by the rows $i_1,i_2,i_3$ of $J$.

Consider the coordinate of $\mathbb{P}_{\K}^{\binom{6}{3}-1}$ corresponding to the position of $\Delta_{1,4,6}$. Dividing the other coordinates by this minor we obtain an affine chart of $X(S)^*$. This affine chart corresponds to the 3-dimensional toric variety defined by $S_A=\mathbb{Z}_{\geq0}(g_1,\ldots,g_8)$. Indeed, this can be verified through straightforward computations: it is just a matter of computing the 20 minors of the Jacobian matrix (all of them monomials) and dividing them by $\Delta_{1,4,6}$; these quotients being monomials, their exponents generate a semigroup which, in turn, is generated by the set $H_1=\{g_1,\ldots,g_8\}$.

For the second iteration of the Nash blowup, we proceed analogously starting with the toric variety $X(S_A)\subset\K^8$, corresponding to $S_A=\mathbb{Z}_{\geq0}(g_1,\ldots,g_8)$. In this case we look for the affine chart of $\mathbb{P}_{\K}^{\binom{8}{3}-1}$ corresponding to the position of $\Delta_{1,5,7}$. This affine chart turns out to be isomorphic to $X(S)$. Again, this is verified through straightforward computations: after computing the 56 minors of the new Jacobian matrix (all of them monomials) and dividing them by $\Delta_{1,5,7}$, the resulting semigroup can be generated by the set $H_2$. As before, we conclude that $X(T)$ is isomorphic to $X(S)$.
\end{remark}

\begin{remark} \label{rem:equation}
Let $\K[x_1,\dots,x_6]$ be the polynomial ring in $6$ variables. The map  $\K[x_1,\dots,x_6]\to \K[S]$ given by $x_i\mapsto \chi^{h_i}$ is surjective. Its kernel $I$  is generated by 
\begin{align*}
\begin{array}{ccccc}
x_5^2 - x_3 x_6,&  x_1 x_5 - x_2 x_6,&  x_1 x_3 - x_2 x_5,&
x_1^2 x_2 x_4 - x_3^2,&  x_1^3 x_4 - x_3 x_5 \,. 
\end{array}
\end{align*}
The singular locus of $X(S)$ in coordinates $x_1, \dots, x_6$ is the union of the following two-dimensional linear spaces: $\{x_1=x_2=x_3=x_5=0\}$ and $\{x_1=x_3=x_5=x_6=0\}$.
\end{remark}

\section{Remarks on computer experimentation}\label{sec computer}

In \cite{CDLAL} we presented counterexamples to the Nash blowup conjecture of toric varieties of dimension greater or equal than four. As explained therein, the counterexamples were obtained by implementing the normalized Nash blowup and iterating this process. Despite intensive computational research, we were unable to find a three-dimensional counterexample to the normalized Nash conjecture. 

In this work, we implemented the algorithm for the Nash blowup of a not necessarily normal toric variety using SageMath~\cite{sagemath}. We provide a brief account here, while a full discussion on the implementation and computational aspects will be presented in a forthcoming work~\cite{CDLLcomputational}. Since there is no built-in class in SageMath to handle affine semigroups that are not saturated, we represent semigroups via their Hilbert basis $\mathcal{H}$, ordered lexicographically.

Now, let $S$ be a pointed affine semigroup and let $\mathcal{H}$ be its Hilbert basis. The Nash blowup of the toric variety $X(S)$ is computed as follows:
\begin{enumerate}
  \item Compute all subsets $A = \{h_{i_1}, \dots, h_{i_d}\}$ of $\mathcal{H}$ such that $\operatorname{det}(h_{i_1} \cdots h_{i_d}) \neq 0$.
  \item Compute $\mathcal{G}_A$ as defined in \eqref{eq:Ga} and check whether the semigroup it generates is pointed.
  \item Compute a Hilbert basis of the semigroup generated by $\mathcal{G}_A$.
\end{enumerate}
Then, as stated below \eqref{eq:Ga}, a set of covering affine charts of the normalized Nash blowup of $X(S)$ is given by the toric varieties $X(S_A)$ for all sets $A$ identified in the first step, whenever $S_A$ is pointed.

Implementing the first and second steps of this algorithm in SageMath is straightforward since it only involves standard operations. Moreover, there is a built-in function in SageMath for verifying that $S_A$ is pointed through the cone it generates. The third step is more involved and requires to extract Hilbert basis of an affine semigroup from a finite generating set. 

We implemented the third step as follows. Given a generating set $\mathcal{G}$, we aim to find the irreducible elements, that is, the elements that cannot be written as the sum of two other elements. In other words, we need to remove the reducible elements. To achieve this, we first compute the set of all $N$-sums of $\mathcal{G}$ recursively, and eliminate from $\mathcal{G}$ all elements that we encounter during this process. We refer to this as the sieving process. We use $N = 6$, a value optimized through computer experimentation. This quickly eliminates many reducible elements. For the remaining ones, we check individually whether each can be written as a sum of the others, which amounts to a feasibility problem in linear programming.
This step is easy to state, but may be computationally difficult to solve. The complexity of this final step depends on the size of the set of surviving elements of $\mathcal{G}$ after the sieve, which is why it is important to eliminate as many elements as possible in the first stage. 

\medskip

We proceed by iterating this algorithm, discarding affine charts that are already smooth at each step. The affine toric variety $X(S)$ is resolved by successive Nash blowups if and only if this iterative process eventually terminates.

The non-normal affine variety presented in the proof of the theorem appears for the first time in the third iteration of the Nash blowup of the non-normal affine toric variety $X(S)$, where $S \subset M = \mathbb{Z}^3$ is the affine semigroup generated by the columns of the following matrix
$$
\left[\begin{array}{rrrr}
1 & 0 & 0 & 1  \\
0 & 1 & 0 & 1  \\
0 & 0 & 1 & -6 
\end{array}\right]\,.
$$

The non-normal affine toric variety $X$ presented in the proof of the theorem is the only three-dimensional counterexample to the Nash conjecture that we have found with the property that the second iteration contains an affine chart isomorphic to itself. We have not found any toric variety exhibiting this property at the first iteration. On the other hand, we have found four non-isomorphic affine toric varieties with this property at the fourth iteration.

As mentioned at the end of the introduction, we would like to emphasize that we have also conducted extensive computational research on the two open cases of the toric Nash blowup conjecture: the non-normalized Nash blowup of two-dimensional toric varieties, and the normalized Nash blowup of three-dimensional toric varieties. To date, no counterexamples have been found, and it seems unlikely that such a counterexample will be discovered using our current methods.

\bibliographystyle{alpha}
\bibliography{ref}

\newcommand{\etalchar}[1]{$^{#1}$}
\begin{thebibliography}{CDLAL25}

\bibitem[ALP{\etalchar{+}}11]{Ataetal}
Atanas Atanasov, Christopher Lopez, Alexander Perry, Nicholas Proudfoot, and Michael Thaddeus.
\newblock Resolving toric varieties with {N}ash blowups.
\newblock {\em Exp. Math.}, 20(3):288--303, 2011.

\bibitem[CDLAL]{CDLAL}
Federico Castillo, Daniel Duarte, Maximiliano Leyton-\'Alvarez, and Alvaro Liendo.
\newblock Nash blowup fails to resolve singularities in dimensions four and higher.
\newblock {\em To appear in Annals of Mathematics}.
\newblock \href{https://annals.math.princeton.edu/articles/22273}{https://annals.math.princeton.edu/articles/22273}.

\bibitem[CDLAL25]{CDLLcomputational}
Federico Castillo, Daniel Duarte, Maximiliano Leyton-\'Alvarez, and Alvaro Liendo.
\newblock On a computational approach to the nash blowup problem.
\newblock 2025.
\newblock Forthcoming article.

\bibitem[CLS11]{cox2011toric}
David~A. Cox, John~B. Little, and Henry~K. Schenck.
\newblock {\em Toric varieties}, volume 124.
\newblock American Mathematical Soc., 2011.

\bibitem[DDR25]{DDR}
Tha\'is~M. Dalbelo, Daniel Duarte, and Maria Aparecida~Soares Ruas.
\newblock Nash blowups of 2-generic determinantal varieties in positive characteristic.
\newblock {\em Pure Appl. Math. Q.}, 21(4):1557--1575, 2025.

\bibitem[DJNnB24]{DJNB}
Daniel Duarte, Jack Jeffries, and Luis N\'u\~nez Betancourt.
\newblock Nash blowups of toric varieties in prime characteristic.
\newblock {\em Collect. Math.}, 75(3):629--637, 2024.

\bibitem[DS25]{DS}
Daniel Duarte and Jawad Snoussi.
\newblock Nash blowups of normal toric surfaces: the case of one and two segments.
\newblock {\em published online in Beitr. Algebra Geom., https://doi.org/10.1007/s13366-025-00805-x}, 2025.

\bibitem[DT18]{DG}
Daniel Duarte and Daniel~Green Tripp.
\newblock Nash modification on toric curves.
\newblock In {\em Singularities, algebraic geometry, commutative algebra, and related topics}, pages 191--202. Springer, Cham, 2018.

\bibitem[Dua14]{DuarSurf}
Daniel Duarte.
\newblock Nash modification on toric surfaces.
\newblock {\em Rev. R. Acad. Cienc. Exactas F\'{\i}s. Nat. Ser. A Mat. RACSAM}, 108(1):153--171, 2014.

\bibitem[Ful93]{fulton1993introduction}
William Fulton.
\newblock {\em Introduction to toric varieties}.
\newblock Number 131 in Annals of mathematics studies. Princeton university press, 1993.

\bibitem[GiZE09]{EGZ}
S.~M. Guse\u~in Zade and V.~\`Ebeling.
\newblock On the indices of 1-forms on determinantal singularities.
\newblock {\em Tr. Mat. Inst. Steklova}, 267:119--131, 2009.

\bibitem[GM12]{GrMi12}
Dima Grigoriev and Pierre~D. Milman.
\newblock Nash resolution for binomial varieties as {E}uclidean division. {A} priori termination bound, polynomial complexity in essential dimension 2.
\newblock {\em Adv. Math.}, 231(6):3389--3428, 2012.

\bibitem[GPT14]{GoTe14}
Pedro~D. Gonz\'alez~P\'erez and Bernard Teissier.
\newblock Toric geometry and the {S}emple-{N}ash modification.
\newblock {\em Rev. R. Acad. Cienc. Exactas F\'is. Nat. Ser. A Mat. RACSAM}, 108(1):1--48, 2014.

\bibitem[GS77]{GS1}
Gerardo Gonzalez~Sprinberg.
\newblock \'{E}ventails en dimension {$2$} et transform\'{e} de {N}ash.
\newblock {\em Publications du Centre de Math\'{e}matiques de l'E.N.S.}, pages 1--68, 1977.

\bibitem[GS82]{GS2}
Gerardo Gonzalez-Sprinberg.
\newblock R\'{e}solution de {N}ash des points doubles rationnels.
\newblock {\em Ann. Inst. Fourier (Grenoble)}, 32(2):x, 111--178, 1982.

\bibitem[GS09]{GS3}
Gerard Gonzalez-Sprinberg.
\newblock On {N}ash blow-up of orbifolds.
\newblock In {\em Singularities---{N}iigata--{T}oyama 2007}, volume~56 of {\em Adv. Stud. Pure Math.}, pages 133--149. Math. Soc. Japan, Tokyo, 2009.

\bibitem[Hir83]{Hi83}
Heisuke Hironaka.
\newblock On {N}ash blowing-up.
\newblock In {\em Arithmetic and geometry, {V}ol. {II}}, volume~36 of {\em Progr. Math.}, pages 103--111. Birkh\"auser, Boston, MA, 1983.

\bibitem[LJR03]{LJ-R}
Monique Lejeune-Jalabert and Ana~J. Reguera.
\newblock The {D}enef-{L}oeser series for toric surface singularities.
\newblock In {\em Proceedings of the {I}nternational {C}onference on {A}lgebraic {G}eometry and {S}ingularities ({S}panish) ({S}evilla, 2001)}, volume~19, pages 581--612, 2003.

\bibitem[Nob75]{Nob75}
A.~Nobile.
\newblock Some properties of the {N}ash blowing-up.
\newblock {\em Pacific J. Math.}, 60(1):297--305, 1975.

\bibitem[Oda83]{oda1983convex}
Tadao Oda.
\newblock {\em Convex bodies and algebraic geometry: an introduction to the theory of toric varieties}.
\newblock Springer, 1983.

\bibitem[Reb77]{Reb}
Vaho Rebassoo.
\newblock Desingularization properties of the nash blowing-up process, 1977.
\newblock Ph. D. dissertation, University of Washington.

\bibitem[{Sag}24]{sagemath}
{Sage Developers}.
\newblock {\em {S}ageMath, the {S}age {M}athematics {S}oftware {S}ystem ({V}ersion 10.4)}, 2024.
\newblock {\tt https://www.sagemath.org}.

\bibitem[Sem54]{Se54}
J.~G. Semple.
\newblock Some investigations in the geometry of curve and surface elements.
\newblock {\em Proc. London Math. Soc. (3)}, 4:24--49, 1954.

\bibitem[Spi90]{Sp90}
Mark Spivakovsky.
\newblock Sandwiched singularities and desingularization of surfaces by normalized {N}ash transformations.
\newblock {\em Ann. of Math. (2)}, 131(3):411--491, 1990.

\bibitem[Spi20]{Sp20}
Mark Spivakovsky.
\newblock Resolution of singularities: an introduction.
\newblock In {\em Handbook of geometry and topology of singularities. {I}}, pages 183--242. Springer, Cham, [2020] \copyright 2020.

\bibitem[Stu96]{Sturm}
Bernd Sturmfels.
\newblock {\em Gr\"obner bases and convex polytopes}, volume~8 of {\em University Lecture Series}.
\newblock American Mathematical Society, Providence, RI, 1996.

\end{thebibliography}

\end{document}